# Thermodynamic and vortic structures of *real Schur flows*

Jian-Zhou Zhu (朱建州)[a)]
*Su-Cheng Centre for Fundamental and Interdisciplinary Sciences, Gaochun,*
*Nanjing 211316, China and*
*Fluid Institute, Nanchang University, Nanchang 330031, China.*

A two-component-two-dimensional coupled with one-component-three-dimensional (2C2Dcw1C3D) flow may also be called a real Schur flow (RSF), as its velocity gradient is uniformly of real Schur form, the latter being the intrinsic local property of any general flows. The thermodynamic and 'vortic' fine structures of RSF are exposed and, in particular, the complete set of equations governing a (viscous and/or driven) 2C2Dcw1C3D flow are derived. The Lie invariances of the decomposed vorticity 2-forms of RSFs in $d$-dimensional Euclidean space $\mathbb{E}^d$ for any interger $d \geq 3$ are also proven, and many Lie-invariant fine results, such as those of the combinations of the entropic and vortic quantities, including the invariances of the decomposed Ertel potential vorticity (and their multiplications by any interger powers of entropy) 3-forms, then follow.

## I. INTRODUCTION

The *Taylor column*s are ubiquitous in rotating fluid systems (c.f., e.g., Ref. 1 for a historical account). A 'column' indicates the two-dimensional structure, and, when referring to the velocity field in the rotation plane, it means the total 3-space velocity has a gradient matrix *uniformly* of the *real Schur form* (RSF), particularly for the formal Taylor-Proudman limit of compressible rotating flows.[2] On the other hand, any real matrices are *similar* to the Schur form,[3] and, has appeared in studies of *local* flow patterns[4] and dynamics.[5] The flow with the velocity gradient being globally of uniform RSF may be called a *real Schur flow* (also RSF without ambiguity in the context). The RSF has a two-component-two-dimensional coupled with one-component-three-dimensional (2C2Dcw1C3D) structure,[5] thus is distinct from the traditional two-dimensional passive scalar or two-dimensional-three-component (2D3C) flow.[2,6] 2C2Dcw1C3D flow appears more naturally in between the 2D3C and the most general three-component-three-dimensional (3C3D) flows. Note that depending on which velocity gradient matrix between the two of transposed form of each other or the upper or lower (quasi-)triagular matrix is used,[3,5] an alternative RSF may also be considered with no essential difference for our purpose.

Note also that while other matrix decompositions may also be applied to matrixes which correspond to some definite physical phenomena, such as the application of the Jordan decomposition to study the growth of density perturbations and its connection with giant molecular cloud,[7] they are not generic in the sense that the decompositon do not always exist for a real matrix such as our velocity gradient $G$ (relevant systematic analysis is of course interesting and may deserve a separate investigation). Every matrix/tensor decomposition has of course its own value, but, to our best knowledge, real Schur transform appears to be the most general for $G$ which in principle can be 'anything' real, which makes us consider RSF as the (anisotropic) base structure that should be the generic corner stone of any flows, including the isotropic turbulence. For example, a helical RSF has been argued to be 'fastened' by the helicity through the Taylor-Proudman effect and has been proposed to play the role of *chiral base flow* for understanding the helicity effect on the compressibility

—————
a)Electronic mail: jz@sccfis.org

of a turbulence with a *weak equivalence principle*,[8] which however may be only a scratch on the surface. The other properties of RSF can also be useful for understanding general turbulence and deserve further clarification, especially for a *strong equivalence principle*, if indeed possible, like the situation of general relativity as indicated in Ref. 8. Here, we focus on the most fundamental thermodynamic structures needed to maintain the 2C2Dcw1C3D dynamics and 'vortic' (related-to-vortex-dynamics) properties from 2C2Dcw1C3D character of the velocity field. In the latter, unlike the conventional constructions (e.g., Refs. 9–13), a set of Lie invariances will be derived from studying which and how many independent components of the vorticity 2-form are invariant.

Numerical investigations, helpful for the programme of constructing the (helical) turbulence strong equivalence principle,[8] of 3-space RSF requires precise formulations to design effective and efficient algorithms with different methods [including the possibility of appropriate 'gene insertion' at the 'micro-scopic' (molecular dynamics) and/or meso-scopic (kinetic theory) levels to obtain the 'macro-scopic' RSF], and an internal mathematical motivation was to address the concerns: whether the two Helmholtz/frozen-in laws found in Ref. 5 is merely a coincidence, or there is something deeper or more essential underlying it; and, how about the structures of other thermodynamic variables (c.f., most recently Ref. 10 and references therein for relevant connections)? The former questions lead us to considering general barotropic ideal RSF in $d$-dimensional Euclidean[14] space $\mathbb{E}^d$ for $d \geq 3$,[15] while the last question will also be addressed with more general situations, resulting in the thermodynamic (II A), 'vortic' (II B) and the mixed (II C) fine structures of RSF, in particular, the complete set of equations (4a,4b,4c) governing a (viscous and/or driven) 2C2Dcw1C3D flow; also, the Lie invariances of the decomposed vorticity 2-forms of RSFs in $d$-dimensional Euclidean space $\mathbb{E}^d$ for any interger $d \geq 3$ are proven in Sec. II B, which indicates the two Helmholtz theorems of the complementary components of vorticity found recently[5] in 3-space RSF is not coincidental, but underlain by a general decomposition theorem, thus essential; and what's more, many Lie-invariant fine results, such as those of the combinations of the entropic and vortic quantities, including the invariances of the decomposed Ertel potential vorticity (and their multiplications by any interger powers of entropy) 3-forms, then follow (II C and III), indicating finer Cauchy-invariant equations than those for general full flows proven recently.[12]

Yet another motivation is related to previous efforts with novel approaches to the physics of passive scalars,[16] fast-rotation limit, i.e., the generalized Taylor-Proudman theorem[2] in $\mathbb{E}^d$ and to the expected new physics and insights of quantum flows,[17,18] which will be discussed in Sec. III after all relevant results are established.

## II. ANALYSIS

Differential forms are special covariant tensors and the vorticity 2-form is simply the antisymmetric second order one which can also be represented with a matrix. For example, the vorticity exterior derivative of the 1-form $U := \sum_{i=1}^{d} u_i dx_i$, corresponding to the velocity vector $\boldsymbol{u} := \{u_1, u_2, ..., u_d\}$, reads, with the index '$,_i$' $\leftrightarrow$ '$\partial_{x_i}$', for $d = 4$ in co-ordinate form

$$
\begin{aligned}
dU = {} & (u_{2,1} - u_{1,2})dx_1 \wedge dx_2 + (u_{3,1} - u_{1,3})dx_1 \wedge dx_3 + \\
& + (u_{4,1} - u_{1,4})dx_1 \wedge dx_4 + (u_{3,2} - u_{2,3})dx_2 \wedge dx_3 + \\
& + (u_{4,2} - u_{2,4})dx_2 \wedge dx_4 + (u_{4,3} - u_{3,4})dx_3 \wedge dx_4,
\end{aligned}
\tag{1}
$$

whose matrix representation, for the components associated to the bases $dx_i \wedge dx_j$, writes

$$
\begin{pmatrix}
0 & \frac{u_{1,2} - u_{2,1}}{2} & \frac{u_{1,3} - u_{3,1}}{2} & \frac{u_{1,4} - u_{4,1}}{2} \\
\frac{u_{2,1} - u_{1,2}}{2} & 0 & \frac{u_{2,3} - u_{3,2}}{2} & \frac{u_{2,4} - u_{4,2}}{2} \\
\frac{u_{3,1} - u_{1,3}}{2} & \frac{u_{3,2} - u_{2,3}}{2} & 0 & \frac{u_{3,4} - u_{4,3}}{2} \\
\frac{u_{4,1} - u_{1,4}}{2} & \frac{u_{4,2} - u_{2,4}}{2} & \frac{u_{4,3} - u_{3,4}}{2} & 0
\end{pmatrix}.
\tag{2}
$$

The matrix representation of $\nabla \boldsymbol{u}$ in $\mathbb{E}^d$ is

$$G = \begin{pmatrix} u_{1,1} & u_{2,1} & u_{3,1} & u_{4,1} & u_{5,1} & u_{6,1} & \cdots & u_{d,1} \\ u_{1,2} & u_{2,2} & u_{3,2} & u_{4,2} & u_{5,2} & u_{6,2} & \cdots & u_{d,2} \\ u_{1,3} & u_{2,3} & u_{3,3} & u_{4,3} & u_{5,3} & u_{6,3} & \cdots & u_{d,3} \\ u_{1,4} & u_{2,4} & u_{3,4} & u_{4,4} & u_{5,4} & u_{6,4} & \cdots & u_{d,4} \\ u_{1,5} & u_{2,5} & u_{3,5} & u_{4,5} & u_{5,5} & u_{6,5} & \cdots & u_{d,5} \\ u_{1,6} & u_{2,6} & u_{3,6} & u_{4,6} & u_{5,6} & u_{6,6} & \cdots & u_{d,6} \\ \cdots & \cdots & \cdots & & & & & \cdots \\ u_{1,d} & u_{2,d} & u_{3,d} & u_{4,d} & u_{5,d} & u_{6,d} & \cdots & u_{d,d} \end{pmatrix}. \tag{3}$$

The three-dimensional (3D) case corresponds to the left-upper $3 \times 3$ block, the RSF of which is arranged to have its left-upper $2 \times 2$ block corresponding to the two conjugate complex eigenvalues of $G$, thus two vanishing left-lower elements (in blue color) indicating a two-component-two-dimensional coupled with one-component-three-dimensional (2C2Dcw1C3D[5]) flow: when all eigenvalues are real, $u_{1,2}$ also vanishes, which is a stronger condition but which is not generic and does not lead to essentially interesting new results in the current discussions, thus will not be particularly discussed; similarly is for $d > 3$. [The RSF is nonunique, depending on how the order of the eigenvalues or the corresponding coordinates are arranged, which however is not essential.] $G$ can be decomposed into symmetric and anti-symmetric parts, $D = (G + G^T)/2$ and $A = (G - G^T)/2$, the latter, given in Eq. (2) for $d = 4$, may be viewed as a representation of the vorticity 2-form $\Omega = dU$.

Let's start with the equation in $\mathbb{E}^3$ for the RSF with density $\rho$, pressure $p$, velocity $\boldsymbol{u}$ and, for simplicity, the dissipation term $D(\boldsymbol{u}) = \nu \nabla^2 \boldsymbol{u}$ with constant kinetic viscosity $\nu$:[19]

$$\partial_t \rho + \nabla \cdot (\rho \boldsymbol{u}) = 0, \tag{4}$$

$$\partial_t \boldsymbol{u}_h + \boldsymbol{u}_h \cdot \nabla_h \boldsymbol{u}_h = -\rho^{-1} \nabla_h p + \nu \nabla_h^2 \boldsymbol{u}_h, \tag{4b}$$

$$\partial_t u_3 + \boldsymbol{u}_h \cdot \nabla_h u_3 + u_3 u_{3,3} = -\rho^{-1} p_{,3} + \nu \nabla^2 u_3, \tag{4c}$$

where $x_1$ and $x_2$ are the '$_h$orizontal' coordinates and the corresponding $\boldsymbol{u}_h := \{u_1, u_2\}$ is independent of the 'vertical' coordinate $x_3$, i.e.,

$$\boldsymbol{u}_{h,3} \equiv 0. \tag{5}$$

Higher-dimensional case can be similarly formulated. For the barotropic case, we can introduce the specific enthalpy $\Pi$ (up to a constant)

$$\nabla \Pi := (\nabla p)/\rho, \tag{6}$$

and the isothermal (constant-temperature) relation

$$p = c^2 \rho \text{ results in } \nabla \Pi = c^2 \nabla \ln \rho, \tag{7}$$

where $c$ is the sound speed.

In Ref. 5, the author noticed that the incompressible RSF admits no periodic solution, but richer structures can present in the compressible one. Indeed, we see no reason to exclude periodic solutions for compressible RSF and actually preliminary numerical tests have shown the realizability of RSF turbulence in a cyclic box, but specific analysis of the latter, especially the numerical investigation of its detailed statistical dynamics, however belongs to another communication for different purpose.[20] Note that the selective cylindrical condition for the horizontal flow also presents in the Taylor-Proudman limit of fast rotating compressible flows,[2] as remarked earlier, thus such a 2C2Dcw1C3D flow is quite of physical sense, instead of being purely artificial. However, the Taylor-Proudman limit has additionally $\nabla_h \cdot \boldsymbol{u}_h = 0$ which is not required by our general RSF.

## A. Thermodynamic structures

The RSF has not only the defining characteristics in the velocity but also some particular thermodynamic structures which governs the dynamics. Below, we will derive the results from the 2C2Dcw1C3D velocity field for the 3-space dynamics.

Taking derivative with respect to $x_3$ in Eq. (4b) for the horizontal momentum, Eq. (5) requires[21]

$$[(\nabla_h p)/\rho]_{,3} = 0 \tag{8}$$

by the requirement of RSF uniformly over space and time, which leads to nontrivial consequence as follows.

### 1. Barotropic structures

By Eq. (6), the barotropic Eq. (8) writes

$$(\nabla_h \Pi)_{,3} = 0 \ (\text{or } \Pi_{,13} = \Pi_{,23} = 0), \tag{9}$$

which simply means, in words, $\Pi$ should be decomposed into two functions, $\mathscr{P}_h$ and $\mathscr{P}_3$, one of only the horizontal coordinate $\boldsymbol{x}_h := \{x_1, x_2\}$ and the other of only $x_3$. To be precise, Eq. (9) means $\boldsymbol{F}_h(\boldsymbol{x}_h) := \nabla_h \Pi$ is a (vector) function of only the horizontal coordinate $\boldsymbol{x}_h$ and $F_3(x_3) := \Pi_{,3}$ a function of only $x_3$, which leads to

$$\Pi = \mathscr{P}_3(x_3) + \mathscr{P}_h(x_1, x_2) \tag{10}$$

where

$$\mathscr{P}_h := \int\int \boldsymbol{F}_h \cdot d^2 \boldsymbol{x}_h \Big|^h \text{ and } \mathscr{P}_3 := \int F_3 dx_3 \Big|^3, \tag{11}$$

with $\Big|^h$ meaning "ignoring the $x_3$-dependent integration constant which is absorbed by $\mathscr{P}_3$" and $\Big|^3$ "ignoring the $\boldsymbol{x}_h$-dependent integration constant which is absorbed by $\mathscr{P}_h$".

We may consider the RSF in a box of dimension $L_z \times L_2 \times L_3$, cyclic in each direction, or with $L \to \infty$ in some direction(s) with the field vanishing sufficiently fast. Introducing

$$\langle \bullet \rangle_{12} := \frac{\int_0^{L_2} \int_0^{L_1} \bullet dx_1 dx_2}{L_1 L_2}, \ \langle \bullet \rangle_3 := \frac{\int_0^{L_3} \bullet dx_3}{L_3} \text{ and } \langle \bullet \rangle_{123} := \frac{\int_0^{L_1} \int_0^{L_2} \int_0^{L_3} \bullet dx_1 dx_2 dx_3}{L_1 L_2 L_3}, \tag{12}$$

we have

$$\langle \Pi \rangle_3 = \mathscr{P}_h(x_1, x_2) + \langle \mathscr{P}_3 \rangle_3, \tag{13a}$$

$$\langle \Pi \rangle_{12} = \mathscr{P}_3(x_3) + \langle \mathscr{P}_h \rangle_{12}, \tag{13b}$$

$$\langle \Pi \rangle_{123} = \langle \mathscr{P}_h \rangle_{12} + \langle \mathscr{P}_3 \rangle_3 = \langle \Pi \rangle_{12} + \langle \Pi \rangle_3 - \Pi,$$

$$\text{i.e., } \Pi = \langle \Pi \rangle_{12} + \langle \Pi \rangle_3 - \langle \Pi \rangle_{123}. \tag{13c}$$

Rewriting Eq. (4) as [c.f., equation (2) of Ref. 8 and references therein]

$$\partial_t \ln \rho = -\boldsymbol{u} \cdot \nabla \ln \rho - \nabla \cdot \boldsymbol{u} \tag{14}$$

and taking, with Eq. (7) for the isothermal case, $\Pi$ in Eq. (10) particularly

$$\Pi = c^2 \ln \rho = c^2(\langle \ln \rho \rangle_{12} + \langle \ln \rho \rangle_3 - \langle \ln \rho \rangle_{123}), \tag{15}$$

and, by bringing Eqs. (15 and 13c) into Eq. (14) with the interchangeability of the order of the time derivative and the average operators defined in Eq. (12) and with

$$\langle u_{1,1} \rangle_1 = \langle u_{2,2} \rangle_2 = \langle u_{3,3} \rangle_3 = 0 \tag{16}$$

from the periodic boundary condition, we finally obtain the partial-integral-differential equation (not the original conventional partial differential equation — PDE)

$$\partial_t \ln \rho = \langle \boldsymbol{u} \cdot \nabla \ln \rho \rangle_{123} - \langle \boldsymbol{u} \cdot \nabla \ln \rho + u_{3,3} \rangle_{12} - \langle \boldsymbol{u} \cdot \nabla \ln \rho + u_{1,1} + u_{2,2} \rangle_3. \tag{4a}$$

It is seen that the above derivation works also for the case with a 2C2Dcw1C3D external force, actually the acceleration of the velocity, $\boldsymbol{a}$, applied to Eqs. (4b and/or 4c), since the horizontal component $\boldsymbol{a}_h$ of such force is independent of $x_3$, keeping still Eq. (8).

Eq. (4a) summarizes Eqs. (4,5 and 8) and, together with Eqs. (4b and 4c), completely defines the dynamics of the 2C2Dcw1C3D isothermal flow in a cyclic box, even for the driven case with 2C2Dcw1C3D force as just mentioned. Actually, we will see that Eq. (4a) applies in more general nonbarotropic RSFs as well. We remark for both applications and further mathematical analysis of equations (say, proving particular properties of the solutions) that such analytical results may be implemented directly in the numerical algorithm or should be realized approximately with desired accuracy at each 'physical' time step with appropriate techniques, say by the standard implicit algorithm with pseudo-time (inner) iterations to converge to the precise relations, the latter being much more expensive but of more general applicability; also possible are the micro- and meso-scopic methods,[22] with specific rules implanted for the desired macroscopic RSF dynamics.

### 2. Nonbarotropic structures

For an ideal gas with, say, $p = \rho \mathcal{R} T$ with $\mathcal{R}$ being a constant, Eq. (8) reads

$$\left[ \frac{\nabla_h(\rho T)}{\rho} \right]_{,3} = 0. \tag{17}$$

Eq. (17) indicates that $\frac{\nabla_h(\rho T)}{\rho}$ is a function of only $x_1$ and $x_2$ and should have such separation of the variables

$$\rho = r(x_3)/R(x_1, x_2) \tag{18}$$

that the numerator and denominator can cancel the common factor $r(x_3)$. In other words, we have the same structure, especially Eq. (4a), as in the isothermal case. And, by taking (18) into (17), we further have

$$T(x_1, x_2, x_3) = \mathcal{T}(x_1, x_2) + R(x_1, x_2)\tau(x_3), \tag{19}$$

where $\mathcal{T}$ and $\tau$ and other variables are of course also time dependent (cf. the first paragraph of Sec. II A 1 for a more precise derivation). Eq. (19) characterizes the fine structures of $T$, which may not appear to be obvious due to simultaneous appearance of $x_1$, $x_2$ and $x_3$ in the second term of the right-hand side (RHS) of Eq. (19), thus we make some elaborations in the following.

From Eq. (19), we have for four (superscript) values of $x_3$

$$\frac{T(x_1, x_2, x_3^2) - T(x_1, x_2, x_3^1)}{T(x_1, x_2, x_3^4) - T(x_1, x_2, x_3^3)} = \frac{\tau(x_3^2) - \tau(x_3^1)}{\tau(x_3^4) - \tau(x_3^3)}; \tag{20}$$

and, for two sets of (superscript) values, respectively, of $x_1, x_2$ and of $x_3$, we have

$$\frac{T(x_1^1, x_2^1, x_3^2) - T(x_1^1, x_2^1, x_3^1)}{T(x_1^2, x_2^2, x_3^2) - T(x_1^2, x_2^2, x_3^1)} = \frac{R(x_1^1, x_2^1)}{R(x_1^2, x_2^2)}. \tag{21}$$

Now the RHS of Eq. (20) is independent of $x_1$ and $x_2$ and the RHS of Eq. (21) is independent of $x_3$, both presenting clear fundamental physical properties that may be useful in mathematical analysis of the equation and numerical simulation and diagnosis.

### 3. Adiabatic Lie-invariant structures

In Ref. 5, the author found an interesting 'vortic' property which is that the frozen-in/Helmholtz theorem of vorticity is decomposed into two invariant laws for the 'horizontal' and 'vertical' vorticities. Logically speaking, it was not clear whether such a frozen-in decomposition was coincidental, happening due to other implicitly used condition(s) such as the dimension number, or essential. Since the general material/frozen-in invariance is described by the Lie invariance, particularly that for the barotropic ideal vorticity $(\partial_t + L_u)dU = 0$, which applies in any dimension number,[11] we thus are curious whether it is underlain by some more essential Lie-invariant theorem for flows in general $d$-dimensional space.

In particular, for adiabatic ideal flow, we have the 0-form/scalar entropy $S$ which is Lie-invariant

$$(\partial_t + L_u)S = 0. \tag{22}$$

Thus, $S^s dU$ with any integer $s$ is still a Lie-invariant 2-form and may be called the "'s-entropic' or 'sentropic' vorticity" (or "vortic 'sentropy'") and the frozen-in decomposition of vorticity in Ref. 5 can be extended to such 2-forms. Similarly, the Ertel 'potential vorticity' 3-form $dS \wedge dU$[10,12,13,23,24] and the "'sentropic' potential vorticity" $S^s dS \wedge dU$ are also accordingly decomposable. But, again, these entropy relevant results, already involving the vorticity, are so far only for flows in $\mathbb{E}^3$, and the essentiality should be examined in $\mathbb{E}^d$ for general $d > 3$, which is a target of the following sections, particularly Sec. II B. Before proceeding, with the already given introductory discussion on the connection between RSF and rotating flows from the beginning, it may be more useful for some readers to further remark that the (exterior) derivative of entropy is associated with the entropy gradient and that may well be related to baroclinic instability in meteorology and relevant disciplinaries:[25,26] our derivation assumes barotropicity (for the Lie invariance of vorticity 2-form) which is sufficient but not always necessary, just as the 'ideal' condition,[27] since the effects of baroclinic ingredients could somehow 'locally' cancel among themselves for $dU$ (but not necessarily for others).

### B. Vortic structures

In terms of differential forms, the inviscid [$\nu = 0$ in Eqs. (4b,4c), say] equations of the horizontal and total velocities read,

$$\partial_t U_h + L_{u_h} U_h = -d_h(\Pi - u_h^2/2), \tag{23}$$

$$\partial_t U + L_u U = -d(\Pi - u^2/2), \tag{24}$$

which also applies for general *ideal* real Schur flows in $\mathbb{E}^d$ with $d \geq 3$, the precise meaning of $U_h$ for $d > 3$ to be further clarified below. Thus, with the interchanges of the exterior derivative with the (partial) time derivative and the Lie derivative, and, with the replacement of $L_{u_h}$ by $L_u$ (Lemma 1), we have the Lie-invariance laws for the vorticity 2-forms $\Omega_h = dU_h$ and $\Omega = dU$,

$$\partial_t \Omega_h + L_u \Omega_h = 0, \tag{25}$$

$$\partial_t \Omega + L_u \Omega = 0. \tag{26}$$

We already can see that the two frozen-in laws of decomposed vorticities (for which the spatial derivatives of $u_h$ and $u_z$ are respectively responsible), found in Ref. 5 for $d = 3$, correspond to Eq. (25) and, according to the linearity of Lie derivative operator, the substraction by it of Eq. (26), but now also for compressible barotropic RSF. [Such a result however does not guarantee a direct extension to $d > 3$, because, as said, in the latter case the precise meaning of the index $_h$ needs special clarification.] The objective of frozen-in decomposition is then translated to Lie-invariant decomposition which can be formulated for general $d$ as the following:

**Definition 1** *A Lie-invariant decomposition of a (Lie-invariant) vorticity 2-form $\Omega$ into $M \geq 1$ components of a barotropic flow in $\mathbb{E}^d$ is that*

$$\Omega = \sum_{i=1}^{M} \Omega_i, \; with \; (\partial_t + L_{\boldsymbol{u}})\Omega_i = 0, \tag{27}$$

*where $\Omega_i$s are linearly independent.*

**Remark 1** *Obviously, $M \leq d(d-1)/2$ (no larger than the number of the elements of the $A$ matrix), and when $M = 1$, $\Omega_1 = \Omega$. Actually, we will see that $M \geq [\frac{d+1}{2}]$, where $[...]$ denotes the integer part.*

### 1. Observation

As remarked in Sec. II B, concerning the left-upper $3 \times 3$ RSF block (with $u_{1,2} = u_{1,3} = 0$, in blue color) designated with single underline and single right 'wall' for $\mathbb{E}^3$ in Eq. (3), the vorticity 2-form, corresponding to the anti-symmetric part of the left-upper $2 \times 2$ block, of the horizontal velocity $\boldsymbol{u}_h$ is Lie-carried by $\boldsymbol{u}_h$, and thus also by the whole 3-space $\boldsymbol{u}$ (Lemma 1). Since the whole $\Omega = d U$ is Lie invariant respect to $\boldsymbol{u}$, the vorticity 2-form component corresponding to anti-symmetric part of the right column of the matrix is also accordingly Lie invariant, from simple subtraction.

Now in the 4-space with the left-upper $4 \times 4$ RSF block, designated with double underlines and double right walls, with extrally $u_{1,4} = u_{2,4} = 0$ (in blue color), the vorticity 2-form component corresponding to the anti-symmetric part of the left-upper $2 \times 2$ block is Lie-carried by the 4-space $\boldsymbol{u}$ (twice applications of Lemma 1). Thus, the vorticity 2-form component corresponding to the anti-symmetric part of the rest two columns on the right, but not the third column as in 3-space, is also accordingly Lie invariant, again, from simple substraction.

Then, in the 5-space, similar to the 3-space case, with the left-upper $5 \times 5$ RSF block with extrally $u_{1,5} = u_{2,5} = u_{3,5} = u_{4,5} = 0$. We see both the vorticity 2-form components, corresponding respectively to the anti-symmetric parts of the left-upper $2 \times 2$ block and to the anti-symmetric part of the third and fourth columns, and their sum as a component of the 5-space vorticity 2-form, are all Lie-carried by $\boldsymbol{u}$ in $\mathbb{E}^5$; thus, again from simple substraction, the vorticity 2-form component corresponding to the anti-symmetric part of the fifth column is also Lie invariant.

The above analysis implies an inductive procedure with 'trivial' (see below for preciseness) extension of dimension and substraction, by which we see how the 6- and higher-space RSF flows carry linearly independent Lie-invariant vorticity 2-form components:

**Theorem 1** *(Lie-invariant vorticity decomposition). With $d+1$ understood to be adding an extra spatial dimension the velocity component of which is constant, i.e., appending an extra column and row of zeros to the bottom and right of the $d$-space velocity gradient matrix $G$, there are $M = [\frac{d+1}{2}]$ linearly independent Lie-invariant (with respect to the ideal barotropic flow in $\mathbb{E}^d$) vorticity 2-form components, each of which subsequently corresponding to the anti-symmetric part of the two columns of $G$ associated to $\Omega_i = dU_i$ for*

$$U_i := u_{2i-1} dx_{2i-1} + u_{2i} dx_{2i}. \tag{28}$$

**Remark 2** *When $i \neq d/2$, each $U_i$ and $\Omega_i$ are respectively perpendicular and parallel to the $d$th coordinate, thus the notations $U_{h_i}$ and $\Omega_{h_i}$, corresponding to Eqs. (23,25), would also be justified.*

### 2. Proof of Theorem 1

**Lemma 1** *The (component of) vorticity 2-form Lie-invariant with respect to a $k$-space velocity is also invariant when the space is trivially extended to $k + 1$ dimensions, where*

*'trivially' refers to the property that the velocity components responsible for the vorticity 2-form do not depend on the extended spatial coordinate.*

With $\boldsymbol{u}$ extended "trivially" from a $k$-space velocity to a $(k+1)$-space $\boldsymbol{u}'$, the proof of Lemma 1 is straightforward by the result of $L_{\boldsymbol{u}} dU_{(i)} = L_{\boldsymbol{u}'} dU'_{(i)}$ derived from $L_{\boldsymbol{v}} = \iota_{\boldsymbol{v}} d + d\iota_{\boldsymbol{v}}$, where $\iota_{\boldsymbol{v}}$ is the interior product with the vector $\boldsymbol{v}$. Thus, we are ready to prove Theorem 1:

*Proof.* With the application of Lemma 1, the induction-fashion argument made in Sec. II B 1 can be adapted with simple replacements of the dimension numbers, 3, 4 and 5 there, with $k = 2n - 1$, $2n$ and $2n + 1$ for $n \geq 2$, to constitute our proof.

Finally, the linear independence of the $\left[\frac{d+1}{2}\right]$ components in Eq. (27) is obvious by the fact that the numbers of bases, $dx_m \wedge dx_n$, involved in different $\Omega_i$s are different, thus the complete proof follows. $\qquad\square$

**Remark 3** *Any linear combinations of $\Omega_i$s are also Lie-invariant. It is very tempting to say that $\Omega_i$s constitute the linearly independent Lie-invariant basis of the invariant vorticity (component) of RSF, but there can be more vanishing components than those left-lower corners of RSF $G$ in Eq. (3), which then may lead to some very special Lie-invariant vorticity component(s) not representable by those $\Omega_i$s: that's why we said in Remark 1 that $M$ could be larger than $\left[\frac{d+1}{2}\right]$.*

### C. Structures of united thermodynamics and vortics

Leaving the remarks concerning those in Sec. II A 3 for higher-order 'mixed' forms, such as the Ertel potential vorticity, to Sec. III and staying still within second-order forms, we have, with $\Omega_i$ given in Theorem 1 and $s$ introduced in Sec. II A 3, the following:

**Corollary 1** *The 'sentropic' vorticity 2-forms, $S^s \Omega_i$, are Lie invariant in the barotropic ideal flows.*

### III. FURTHER DISCUSSIONS

The conventional construction of higher-order Lie invariants then can be applied to the decomposed objects; that is, we have, with $U_i$ and $\Omega_i = dU_i$ in Theorem 1, the following:

**Corollary 2** *Any linear-sum and wedge-product combinations of $\Omega_i$s, and, particularly "'sentropic' potential vorticities" defined by the wedge products of $S$, $dS$ and $\Omega_i$, say, $S^s dS^m \wedge (\Omega_i + \Omega_j)^n$, are also locally invariant.*

Other higher-order pure thermodynamic quantities (Ref. 10 and references therein) accordingly have delicate structures, which however appears not demonstrable in an explicit way so far and is left for future study.

Various new Cauchy invariants equations associated to the decomposed components can be further established by further applying the main theorem of Besse & Frisch,[12] which will explicitly present the fine Lagrangian structures. As we remarked on the thermodynamic fine structures, such fine results may also be useful for Lagrangian numerical techniques and for check of the preciseness of the simulations.

Besides the multi-disciplinary background of studying high-dimensional flows, such as the classical and quantum turbulence, mathematical analysis (of PDEs) and laws of Nature indicated in Ref. 15, it is not impossible that real Schur flows in $\mathbb{E}^d$ with $d > 3$ may be of relevance to applications of fluid dynamics in laboratory physics: advected 2D or 3D passive scalar(s) may be considered as the $d - 2$ or $d - 3$ velocity component(s) independent on the corresponding coordinate(s), i.e., with 'cylinder condition'[16] which however does not simply result from the fast-rotation limit of $d$-space flows,[2] and we are motivated to find extra possible mechanism(s) to reduce from RSF (or its variants with similarly 'nice' properties) to such passive scalar(s). Indeed, the general 'hulking' $d$-space 'full' flow is not 'intimate' or 'sharp' enough, and we did not benefit very much from it in Ref. 16, with a purpose of lifting the dimension to $d$ of a flow to study the scalar(s) advected by this flow by benefitting from the knowlege of $d$-space flow with a 'wholistic' treatment. Now, the RSF or its appropriate variant, just as the 2C2Dcw1C3D flow lies inbetween the 2D3C and 3C3D ones, is not only 'closer' to the passive-scalar problem but also enriched with our finer, thus presumably more powerful, results to be exploited for sharper physical insights.

Using matrix decomposition or transformation techniques, one may perform various other decompositions of the flow and/or the vorticity, which, if being lack of the consideration on the dynamics (such as the Lie invariance with respect to the flow discussed in the above), may be called 'kinetic'. For example, besides the well-known decomposition of $G$ into the symmetric and antisymmetric parts (respectively, $D$ and $A$), Ref. 4 shows that the RSF $G$ can be further transformed into the canonical form and then decomposed into the shear part $S$ and the canonical part $N$, the latter further being composed of the dilation part $E$, the part for the strain rate along some eigenvector $Z$ and the rotation part $\Psi$. For the antisymmetric $A$ representing the vorticity, its canonical form is (block) diagonal with some $2 \times 2$ antisymmetric block(s) and all other elements vanishing:[28] each $2 \times 2$ antisymmetric block represents the (rate of) rotation in the plane (c.f., equation 23 in Ref. 2 and the discussion following it), and, if more than one, all the rotation planes are orthogonal to each other. We call such rotations 'pure'. For example, if the left-lower and right-upper $2 \times 2$ blocks of the matrix (2) are zeros, then the (pure) rotations are respectively in the $x_1$-$x_2$ and $x_3$-$x_4$ planes in the corresponding canonical coordinate frame. However, the $A$ of RSF $G$ is in general not in the canonical form, thus containing 'entangled' rotations. For example, $G$ being of RSF though, none of the nondiagonal elements of matrix (2) is indicated to be vanishing, thus the possibility of simultaneous rotations in $x_1$-$x_2$, $x_3$-$x_4$, $x_1$-$x_3$ and $x_2$-$x_4$ planes in the corresponding coordinate frame. The entanglement of the rotations may lead to ambiguity and even confusions about vorticity, as already presents in 3D case: more than one rotation planes are in general involved in the vorticity, thus very complicated flow pattern from the latter, unless in the canonical frame where there always is only one plane for the pure rotation.

We conclude by returning to remark on the 'dynamical' decomposition: the study of compressibility reduction of helical turbulence in Ref. 8 has made use of helical RSF as the chiral base flow but has not yet exploited the ideal vorticity frozen-in decomposition, however we expect that the latter may be intrinsic to the fundamental mechanisms of other issues of *fully $d$*-dimensional flows (not RSF), including 3D incompressible turbulence. And, carrying such a paradigm over to quantum flows, say, that of the quantum Navier-Stokes model,[17] should be particularly interesting, concerning the relations with the rotating superfluid or Bose-Einstein condensate which may be described by the Gross-Pitaevskii equation (or the 'inviscid' quantum Navier-Stokes, after the Madelung transform),[18] for instance: the quantum pressure term and viscosity may lead to nontrivial differences in the corresponding RSFs, which deserves further investigation.

## ACKNOWLEDGMENTS

Partially supported by NSFC (Nos. 11672102). Discussions with X. He, H. Ren and C.-X. Tang on RSF simulation have been very motivative to establish some of the results. The referee-report comments helped to improve the presentation.

## DATA AVAILABILITY

The data that support the findings of this study are available from the corresponding author upon request.